\input amstex
\UseAMSsymbols
\input epsf.tex
\documentstyle{amsppt}


\loadbold \nologo \pageheight{8.5truein} \pagewidth{7.0truein}
\topmatter
\title Correction to "$C_1$ in [2] is zero"
\endtitle
\author Abbas Bahri
\endauthor
\endtopmatter
John Morgan and G,Tian pointed out a mistake in the concluding
argument for our paper entitled "$C_1$ in [2] is zero", which was
recently published in arXiv:1512.02098. We hereby acknowledge this
mistake and correct the computation, leading to the conclusion that
$C_1$ is non-zero and that their reference [2] does indeed fully
address and resolve the counter-example which we provided in [3] to
the inequality (19.10) in their monograph [1].

\bigskip

The mistake takes place when computing $(\hat{\nabla}_{\frac
{\partial}{\partial t}}\hat{\nabla}_SS,H)$. The metric is variable
here and the derivatives of the Christoffel symbols lead to a
non-zero $C_1$. Namely:

 $(\hat{\nabla}_{\frac {\partial}{\partial
t}}\hat{\nabla}_SS,H)=(\hat{\nabla}_{\hat{H}}\hat{\nabla}_SS,H)-(\hat
{\nabla}_H\hat{\nabla}_SS,H)$

Now, $\hat {H}$ is along $(c(x,t),t)$. Thus, the metric is $g(t)$
and $\hat{\nabla}_SS=H+Ric(S,S)\frac {\partial}{\partial t}$. Thus,

$(\hat{\nabla}_{\hat{H}}\hat{\nabla}_SS,H)=(\hat{\nabla}_{\hat{H}}(H+Ric(S,S)\frac
{\partial}{\partial t}),H)=(\hat{\nabla}_{\frac {\partial}{\partial
t}}(H+Ric(S,S)\frac {\partial}{\partial t}),H)+(\hat
{\nabla}_H(H+Ric(S,S)\frac {\partial}{\partial t}),H)$

Since, $\hat {\nabla}_{\frac {\partial}{\partial t}}\frac
{\partial}{\partial t}=0$ and $(\frac {\partial}{\partial t},H)=0$,
since $(\hat {\nabla}_H(Ric(S,S)\frac {\partial}{\partial
t}),H)=O(k^2)$,

we find that:

$$(\hat{\nabla}_{\frac {\partial}{\partial
t}}\hat{\nabla}_SS,H)=(\hat{\nabla}_{\frac {\partial}{\partial
t}}H,H)+(\hat {\nabla}_HH,H)-(\hat
{\nabla}_H\hat{\nabla}_SS,H)+O(k^2)=$$
$$=(\hat {\nabla}_HH,H)-(\hat
{\nabla}_H\hat{\nabla}_SS,H)+O(k^2)$$

Now, since $S$ is horizontal, $\hat{\nabla}_SS=\nabla_SS+\theta
\frac {\partial}{\partial t}$, $\theta$ bounded, so that

$$(\hat
{\nabla}_H\hat{\nabla}_SS,H)=\hat {\nabla}_H\nabla_SS,H)+O(k^2)$$

Thus, our above expression is, up to $O(k^2)$:

$$(\nabla_HH,H)-(
\nabla_H\nabla_SS,H)$$

$H(c(x,t),s)$ is equal to $\nabla_S^{g(t)}S$, with
$S(c(x,t),s)=\frac {\frac {\partial c(x,t)}{\partial x}}{|\frac
{\partial c(x,t)}{\partial x}|_{g(t)}}$. Along $H$, $(c(x,t),s)$
changes after the time $\tau$ into $(c(x,t+\tau),s)$. With $s=t$,
the metric is $g(t)$, so that, along a piece of curve tangent to $H$
as defined here:

$$\nabla_SS(c(x,t+\tau),s)=\nabla^{g(t)}_SS$$

, with $S(c(x,t+\tau),t)=\frac {\frac {\partial
c(x,t+\tau)}{\partial x}}{|\frac {\partial c(x,t+\tau)}{\partial
x}|_{g(t+\tau)}}$ instead of
$\nabla_SS(c(x,t+\tau),s)=\nabla^{g(t+\tau)}_SS$,with
$S(c(x,t+\tau),t)$ as above. This is the expression that we would
find in $(\nabla_HH,H)$ and there is therefore a difference between
$H(c(x,t+\tau),t)$ and $\nabla^{g(t}_SS$, where $S$ is taken at
$(c(x,t+\tau),t)$. The difference appears through the Christoffel
symbols of the two different metrics $g(t+\tau)$ and $g(t)$. In
$(\nabla_HH,H)-(\nabla_H\nabla_SS,H)$, this difference is
differentiated along $H$, that is along $\tau$ and it leaves a
single factor for $H$, giving rise to $C_1k$, with $C_1$ non-zero.

The observations of John Morgan and Gang Tian, leading to the
complete resolution of this matter, are gratefully acknowledged
here.
\newpage

\widestnumber\no{99999}

\font\tt cmr12 at 24 truept 

\end{document}